\documentstyle[12pt]{article}
                                                 \newif\iffigs\figsfalse

\figstrue

\iffigs
\input epsf
\fi

 \stop \fi

\batchmode
  \newfont{\footbbbfont}{msbm10}
\errorstopmode

\newif\ifamsf\amsftrue
\ifx\footbbbfont\nullfont
  \amsffalse
\fi

\def\ppnumber{duk-m-94-05}
\def\ppdate{April, 1995}

\def\pplogo{\vbox{\kern-\headheight\kern -5pt
\halign{##&##\hfil\cr&{\sc \ppnumber}\cr\rule{0pt}{2.5ex}&\ppdate\cr}
}}
\makeatletter
\advance\textheight by -\headheight
\advance\textheight by -\headsep
\newcommand{\text}[1]{\mathchoice{\mbox{\rm #1}}{\mbox{\rm #1}}
    {\mbox{\scriptsize\rm #1}}{\mbox{\tiny\rm #1}}}

\newif\iffn\fnfalse

\ifamsf
  \newfont{\bigbbbfont}{msbm10 scaled\magstep2}
  \newfont{\bbbfont}{msbm10 scaled\magstep1}  
  \newfont{\smallbbbfont}{msbm8}
  \newfont{\tinybbbfont}{msbm6}
  \newfont{\smallfootbbbfont}{msbm7}
  \newfont{\tinyfootbbbfont}{msbm5}
  \newcommand{\Bbb}[1]{\iffn
      \mathchoice{\mbox{\footbbbfont #1}}{\mbox{\footbbbfont #1}}
      {\mbox{\smallfootbbbfont #1}}{\mbox{\tinyfootbbbfont #1}}\else
      \mathchoice{\mbox{\bbbfont #1}}{\mbox{\bbbfont #1}}
      {\mbox{\smallbbbfont #1}}{\mbox{\tinybbbfont #1}}\fi}
\else
  \def\bigbbbfont{\bf}
  \def\Bbb{\bf}
\fi

\newcommand{\operatorname}[1]{\mathop{\rm #1}\nolimits}

\newcommand{\pfit}[1]{{\it #1:\/}}
\newcommand{\qedsymbol}{Q.E.D.}
\newcommand{\eqref}[1]{(\ref{#1})}

\newenvironment{pf*}{\noindent\pfit }{\qedsymbol\par\par\medskip\par}

\newenvironment{bmatrix}{\left[\begin{array}{ccccc}}{\end{array}\right]}

\def\title{\@dblarg\@xtitle}
\def\@title{}
\def\@xtitle[#1]#2{\gdef\@title{#2}}
\def\author{\@dblarg\@xauthor}
\def\@author{}
\def\@xauthor[#1]#2{\gdef\@author{#2}}
\def\@address{}
\def\address#1{\gdef\@address{#1}}
\def\@email{}
\def\email#1{\gdef\@email{#1}}

\def\trailer{\bigskip\par\noindent \footnotesize {\sc \@address}
   \par\noindent {\it E-mail address:} \@email}
\date{}
\def\dedicatory#1{\def\@date{\normalsize\it#1}}
\def\subjclass#1{\def\@thefnmark{}\@footnotetext{1991
    {\it Mathematics Subject Classification.} #1}}
\def\keywords#1{\def\@thefnmark{}\@footnotetext{
    {\it Key words and phrases.} #1}}
\def\thanks#1{\def\@thefnmark{}\@footnotetext{#1}}
\def\ps@firstpage{\ps@plain \def\@oddhead{\hss\pplogo}%
  \let\@evenhead\@oddhead 
}
\def\maketitle{\par
 \begingroup
 \def\thefootnote{\fnsymbol{footnote}}
 \def\@makefnmark{\hbox
 to 0pt{$^{\@thefnmark}$\hss}}
 \if@twocolumn
 \twocolumn[\@maketitle]
 \else \newpage
 \global\@topnum\z@ \@maketitle \fi\thispagestyle{firstpage}\@thanks
 \endgroup
 \setcounter{footnote}{0}
 \let\maketitle\relax
 \let\@maketitle\relax
 \gdef\@thanks{}\gdef\@author{}\gdef\@title{}\let\thanks\relax}

\makeatother

\newtheorem{definition}{Definition} 

\newcommand{\dlog}{d\mskip0.5mu\log}
\renewcommand{\Im}{\operatorname{Im}}
\newcommand{\Aut}{\operatorname{Aut}}
\newcommand{\C}{{\Bbb C}}
\newcommand{\R}{{\Bbb R}}
\newcommand{\Z}{{\Bbb Z}}
\newcommand{\Hom}{\operatorname{Hom}}
\newcommand{\U}{\operatorname{U}}
\newcommand{\p}{{\Bbb P}}

\newcommand{\Ker}{\operatorname{Ker}}
\newcommand{\SU}{\operatorname{SU}}

\def\lhk{\mathbin{
\hbox{\vrule height1.4pt width4pt depth-1pt
\vrule height4pt width0.4pt depth-1pt}}}

\makeatletter
\relax
\@writefile{toc}{\string\contentsline\space {section}{\string\numberline\space
{1}Coordinates on the B-model moduli space}{3}}
\@writefile{toc}{\string\contentsline\space {section}{\string\numberline\space
{2}The large radius limit}{5}}
\@writefile{toc}{\string\contentsline\space
{subsection}{\string\numberline\space {2.1}The nonlinear $\sigma $-model}{5}}
\@writefile{toc}{\string\contentsline\space
{subsection}{\string\numberline\space {2.2}The A-model parameter space}{6}}
\@writefile{toc}{\string\contentsline\space
{subsection}{\string\numberline\space {2.3}Flat coordinates and the large
radius limit}{7}}
\@writefile{toc}{\string\contentsline\space {section}{\string\numberline\space
{3}Maximally unipotent monodromy}{9}}
\@writefile{toc}{\string\contentsline\space {section}{\string\numberline\space
{4}Equivalence among boundary points}{13}}
\@writefile{toc}{\string\contentsline\space {section}{\string\numberline\space
{5}Determining the mirror map (two conjectures)}{15}}
\@writefile{toc}{\string\contentsline\space
{subsection}{\string\numberline\space {5.1}A conjecture about integral
cohomology}{16}}
\@writefile{lot}{\string\contentsline\space {table}{\string\numberline\space
{1}{\ignorespaces Monodromy calculations}}{18}}
\@writefile{toc}{\string\contentsline\space
{subsection}{\string\numberline\space {5.2}The monomial-divisor mirror
map}{19}}
\@writefile{toc}{\string\contentsline\space {section}{\string\numberline\space
{6}Making enumerative predictions}{19}}
\makeatother

\begin{document}

\title[Making enumerative predictions]{Making enumerative predictions\\
by means of mirror symmetry}

\author
{David R. Morrison}

\thanks{Research partially supported by NSF grant DMS-9401447.}

\address
{Department of Mathematics,
Box 90320,
Duke University,
Durham, NC 27708-0320}

\email
{drm@math.duke.edu}

\date{}

\maketitle

\begin{abstract}
Given two Calabi--Yau threefolds which are believed to constitute
a mirror pair, there are very precise predictions about the
enumerative geometry of rational curves on
one of the manifolds which can be made by
performing calculations on the other.  We review the mechanics
of making these predictions, including a discussion of
two conjectures which specify
how the elusive ``constants of integration'' in the mirror map should be fixed.
Such predictions can be useful for checking whether or not various
conjectural constructions of mirror manifolds are producing reasonable
answers.
\end{abstract}

Two-dimensional quantum field theories constructed from Calabi--Yau
manifolds have been the subject of intensive study over the last
several years.  One of the most intriguing features of these
models is the phenomenon known as
{\em mirror symmetry}\/ \cite{dixon,LVW,CLS,GP}, in
which it is observed that certain pairs of Calabi--Yau manifolds
 produce quantum field theories which appear to be isomorphic via
an automorphism which changes the sign of a certain ``$\U(1)$-charge''
in the theory.
This observation has provided physicists with a powerful computational
tool, since calculations which are difficult in one realization of
the quantum field theory may become much easier in the other (thanks
to a significant shift in geometric interpretation which accompanies
the sign change of the $\U(1)$-charge).
The calculations in question can often be formulated in purely mathematical
terms,
but it should be borne in mind that
the arguments in favor of the equivalence of the answers (when the
calculations are performed on
two different members of a mirror pair)
 rely upon path integral methods which have not yet
been made mathematically rigorous.  For this reason, mathematicians currently
regard these calculations as {\em predicting}\/ rather than
{\em establishing}\/ the results.

The most common application of mirror symmetry as a computational
tool in physics has involved
 calculations of certain so-called ``topological'' correlation
functions of the physical theory.  There are two types of such functions,
associated to the two topological quantum field theories known
 as the ``A-model'' and the ``B-model'' \cite{witten-mirror},
and mirror symmetry predicts that a B-model
calculation on one Calabi--Yau manifold should produce the answer
 to a mirror A-model problem.  Since that answer typically
involves enumerative geometry on the mirror partner, we refer to this
process as ``making enumerative predictions''.

In this paper, we review the mechanics of making these enumerative
predictions (including a discussion of the reasoning which leads to the methods
we describe).  The starting point is a {\em candidate mirror pair}, that
is, a pair of Calabi--Yau manifolds $(X,Y)$ which there is some reason
to believe ought to be a pair of mirror manifolds.
Precise enumerative predictions will be derived from the assumption
that these manifolds are indeed a mirror pair.  If those enumerative
predictions can be verified, the verifications constitute evidence that
a mirror partner has been correctly identified.  Conversely, if
discrepancies are discovered
between enumerative predictions and actual enumerative calculations,
the validity of either a proposed mirror construction or of the precise
geometric
interpretation of the quantum field theories will be called into question.
In the by now rather long list of papers which have made specific
enumerative predictions and attempted to verify them
\cite{CdGP,katzconics,ES,picard-fuchs,font,kt1,lib-teit,kt2,BvS,2param1,%
HKTY1,2param2,HKTY2,chiral,Voisin},
there have been no discrepancies (once all the dust has settled\footnote{I
would
like to take this opportunity to acknowledge an error in \cite{picard-fuchs}
(previously pointed out in \cite{font} and \cite{kt1}):  the normalization
in the last example of that paper was incorrect, and all entries for
coefficients
$n_j$ in the last
row of table 3 of \cite{picard-fuchs} should be divided by 2.}).

We will restrict our attention to the enumerative problem of counting
rational curves on Calabi--Yau threefolds.
The very interesting extensions to higher genus \cite{BCOV,BCOV2}
and to higher dimension \cite{higherD,katz2,kanter}
will not be included.

The general strategy for making enumerative predictions is as follows.
We start with a Calabi-Yau threefold $X$ and a candidate mirror partner $Y$
of $X$.  We formulate enumerative predictions for $X$
by making Hodge-theoretic
calculations on $Y$ (which are a purely mathematical version\footnote{The
connection between the physical and mathematical versions of these
calculations is reviewed in detail in \cite{ceresole}.} of
the B-model calculations mentioned above).  To make these predictions,
 we wish to compare  the asymptotic
expansion of the A-model correlation functions of $X$ with an
appropriate expansion of B-model correlation functions of $Y$.
This comparison is guided by two key observations:
\begin{enumerate}
\item The A-model
comes with distinguished coordinates (the so-called ``flat'' coordinates)
which correspond under mirror symmetry to coordinates on the B-model side
given by ratios of
periods,
and

\item  The coordinates on the A-model side are only well-defined up to
translations by an integral lattice, which implies that the ratios of
periods on the B-model side must exhibit a similar ambiguity
(which will come from the {\em monodromy}\/ of the periods).
\end{enumerate}
Finding ratios of periods whose monodromy behavior (in a particular
open set) reproduces the type of
ambiguity expected from the A-model
then identifies where the expansion of the B-model correlation
function should be made.

We work in this paper
 with completely general Calabi--Yau threefolds, although virtually all
of the specific
calculations to date in the literature have dealt with special
cases (complete intersections in toric varieties)
which are closely related to hypergeometric functions.
We will point out some of the special techniques which are available in those
cases, although we will not review those techniques in detail.

\section{Coordinates on the B-model moduli space} \label{s:bcoords}

A {\em Calabi--Yau threefold}\/ is a compact oriented $6$-manifold $Y$ which
admits
Riemannian metrics whose (global) holonomy is contained in $\SU(3)$.
For any such metric, there exists at least one complex structure
with respect to which the
metric is K\"ahler, and for each such complex structure ${\cal J}$
there is a nowhere-vanishing
holomorphic $3$-form $\Omega$ on the complex manifold $Y_{\cal J}$.

Given a Calabi--Yau threefold, the topological quantum field theory known
as the {\em B-model of $Y$}\/ has as its essential parameters the choice of
complex structure $\cal J$ on $Y$.  In fact we should identify
the {\em B-model moduli space}\/
 with the usual  moduli space of complex structures (with trivial canonical
bundle)
which is studied in
algebraic geometry.  There are some well-known technical difficulties in
constructing such  moduli spaces, but we are primarily
concerned with two aspects of the moduli problem: we need to understand
the moduli space {\em locally}, and we need to be sure that there
are {\em good compactifications}\/ of the moduli space.
For Calabi--Yau manifolds,
the first aspect is covered by the theorem of Bogomolov, Tian and Todorov
\cite{bogomolov,tian,todorov}, which says that the moduli
space is smooth and that its tangent space at $\cal J$ can be naturally
identified with $H^1(T^{(1,0)}_{Y_{\cal J}})$, the first cohomology
group of the sheaf of holomorphic
vector fields.  The second aspect---the existence of a good
compactification---follows from Viehweg's theorem \cite{viehweg}
that the moduli space of {\em polarized}\/ Calabi--Yau manifolds
is a quasi-projective variety.  To compactify the moduli space,
take a projective completion of one of Viehweg's spaces.\footnote{Because
of the polarization condition, the resulting space is only a compactification
of an {\em open subset}\/ of the original moduli space, but this is adequate
for our purposes.}
Other compactifications can then be found by blowing up the original one.

The physically natural coordinates on this B-model moduli space are provided
by ratios of periods of (any) holomorphic $3$-form $\Omega$.  That is,
if $\dim H^1(T^{(1,0)}_{Y_{\cal J}})=r$ then we choose
 $r{+}1$ elements $\gamma_0$, $\gamma_1$, \dots, $\gamma_r$
in $H_3(Y,\Z)$, and  use the ratios
$\int_{\gamma_j} \Omega/\int_{\gamma_0} \Omega$
as local coordinates.  (This form of the coordinates was arrived at empirically
in \cite{CdGP} and explained in terms of conformal field theory in
\cite{BCOV2}.)
At each point in the moduli space, any generic choice of such ratios
will provide good local coordinates, thanks to the local Torelli theorem,
the Bogomolov--Tian--Todorov theorem cited above, and the analysis by
Bryant and Griffiths \cite{bryant-griffiths} of the period map for
such variations of Hodge structure.

One important aspect, therefore, of the problem of making enumerative
predictions
will be to calculate such periods.
This can sometimes be done directly, but a more common approach is an indirect
one in which one first calculates the differential equations which the periods
satisfy.
Choose a family
of holomorphic $3$-forms $\Omega(s)$ which depends on a parameter $s$ on the
moduli space.  (This can only be done locally on the moduli space.)
The periods $\int_{\gamma} \Omega(s)$
can be differentiated with respect to parameters, and there must be
differential
operators ${\cal D}$
 which annihilate the periods, that is
\[
{\cal D}\left(\int_{\gamma} \Omega(s)\right)=0
\]
 for all $\gamma\in H_3(Y,\C)$. (These form a differential ideal on the
moduli space.) We call these differential operators
the {\em Picard--Fuchs operators},
and call the resulting differential equations ${\cal D}\varphi=0$
the {\em Picard--Fuchs equations}\/ determined by $\Omega(s)$.

In principle, the Picard--Fuchs equations are derived as follows.  Let
$\pi:{\cal Y}\to S$ be a proper holomorphic map such that each fiber
$\pi^{-1}(s)$ is a complex manifold diffeomorphic to $Y$ which has
the complex structure corresponding to $s\in S$.  (Such ``universal families''
should
at least exist locally over the moduli space.)  Let
${\Bbb V}=R^3\pi_*\C_{\cal Y}$
be the local system of cohomology groups, and let
${\cal V}={\Bbb V}\otimes {\cal O}_S$ be the corresponding locally free sheaf,
with flat connection
\[\nabla:{\cal V}\to{\cal V}\otimes T^*_S\]
which annihilates sections of ${\Bbb V}$.
This {\em Gauss--Manin connection}\/ can actually be computed in purely
algebraic terms \cite{Katz-Oda}.
Doing so leads to the Picard--Fuchs equations indirectly: if
we choose a basis $\gamma_0, \dots, \gamma_{2r+1}$ of $H_3(Y)$ then we
can write
\[\Omega(s)=\sum\left(\int_{\gamma_j}\Omega(s)\right)e^j,\]
where $e^j$ is the dual basis of cohomology.  Then
\[\nabla \Omega(s)=\sum\left(d\int_{\gamma_j}\Omega(s)\right)e^j.\]
Thus, we can calculate the effect of differential operators on the
periods by calculating the effect of the Gauss--Manin connection
on the cohomology itself, and thereby determine the Picard--Fuchs
equations.

In practice, calculating either the Gauss--Manin connection or the
Picard--Fuchs
equations is rather difficult.  The cases in which these calculations have
been carried out explicitly have involved one of two techniques:
\begin{enumerate}
\item In some cases
 it has been possible to explicitly evaluate some particular
period integral, and expand its value in a power series.  The Picard--Fuchs
equations can then be found by finding which differential operators annihilate
this known period.
This method was
pioneered in \cite{CdGP}, applied in
\cite{CdGP,font,kt1,2param1,2param2,BK}, and reached its culmination
in \cite{periods}.
\item In somewhat greater generality, in many cases
it has been possible to identify the
periods with certain generalized hypergeometric functions; the Picard--Fuchs
equations are then related to the differential equations of
Gel'fand--Zelevinsky--Kapranov \cite{GZK:h}.  This method was first suggested
in
\cite{Bat:vmhs}, developed in \cite{BvS,small}, and systematized
in \cite{HKTY1,HKTY2}.
\end{enumerate}
We refer the interested reader to the cited papers
for more details concerning this part of the calculation.

In general, the Picard--Fuchs equations will have a $(2r{+}2)$-dimensional
family of local solutions at any point of the moduli space, corresponding to
the possible homology classes $\gamma$.  The reduces the problem of
identifying appropriate coordinates to the problem of selecting
the ``correct'' homology classes  $\gamma_0$ and other  $\gamma_j$'s.
We address this problem in
the remainder of this paper.
We will identify the ``correct'' homology classes
 by comparison with the behavior of the A-model, to which we now turn.

\section{The large radius limit} \label{s:largeradius}

The flat coordinates in the A-model---mirror to the ``ratio of periods''
coordinates
discussed in the previous section---have an ambiguity in their definition
which can be described in terms of an integral lattice.
In order to explain this, we first review some of the mathematical
aspects of the moduli spaces
of nonlinear $\sigma$-models (cf.\ \cite{compact,icm}), which
involve both A-model and B-model parameters.

\subsection{The nonlinear $\sigma$-model}

We briefly recall the Lagrangian formulation of nonlinear $\sigma$-models
in dimension 2.  The essential ingredients needed to describe a nonlinear
$\sigma$-model consist of
 a compact manifold $X$, a Riemannian metric $g_{ij}$ on
$X$, and a class $B\in H^2(X,\R/\Z)$,
all defined up to diffeomorphisms
of $X$.
(We represent $B$ as a closed, $\R/\Z$-valued $2$-form, that is, as
a collection of locally defined closed real $2$-forms, the union of whose
domains of
definition is all of $X$, such that the difference between any two local
representatives is
$\Z$-valued wherever it is defined.)
The nonlinear $\sigma$-model is then constructed from
 a $\C/\Z$-valued (Euclidean) action ${\cal S}$
whose bosonic part
assigns to each sufficiently smooth
 map $\phi$
from an oriented Riemannian $2$-manifold $\Sigma$ to $X$
the quantity\footnote{We suppress
the string coupling constant, and use a normalization in which
the action appears as $\exp(2\pi i{\cal S})$
in the partition and correlation functions.}
\[
{\cal S}_{\text{bosonic}}[\phi]:=
i\int_\Sigma \|d\phi\|^2\,d\mu+\int_\Sigma \phi^*(B),
\]
where the norm $\|d\phi\|$ of $d\phi\in\Hom(T_\Sigma,\phi^*(T_X))$
is determined
from the Riemannian metrics on $X$ and on $\Sigma$, and where
$\int_\Sigma \phi^*(B)$ is a well-defined element of $\R/\Z$ by
virtue of the canonical isomorphism $H^2(\Sigma,\R/\Z)\cong\R/\Z$.
(Additional fermionic terms must be added to
${\cal S}_{\text{bosonic}}$
in order to make the theory supersymmetric,
but as they do not affect the essential parameters in the theory we
suppress them here.)

It is more customary to require $B$ to be a
real $2$-form, in which case
${\cal S}_{\text{bosonic}}$
becomes $\C$-valued, and
 one observes that the physics is
invariant under shifting $B$ by an integral cohomology class.
(The possibility of a more general form of the action\footnote{The
more general form of the action and its properties as described in
this paragraph
 arose in discussion with Paul Aspinwall (cf.\ \cite{stable}).}
${\cal S}_{\text{bosonic}}$
which allows $B$ to be an $\R/\Z$-valued $2$-form is
implicit in \cite{Vafa,DijkgraafWitten,chiral}.)
To compare this more
general form to the customary one, consider the exact sequence
\[
0\to H^2_{\text{DR}}(X,\Z)
\to H^2(X,\R) \to H^2(X,\R/\Z)
\to H^3(X,\Z)_{\text{tors}} \to0,
\]
where $H^2_{\text{DR}}(X,\Z)$ denotes the image of $H^2(X,\Z)$ in
de~Rham cohomology.  The last term in
this exact sequence
is a finite group
which labels the connected components of $H^2(X,\R/\Z)$.
If we only used real $2$-forms modulo integral $2$-forms to describe
$B$, we would get only
one connected component of that space.

We are interested in a special case of this construction in which the
theory has what is called $N{=}(2,2)$ supersymmetry and is in addition
conformally invariant.  To ensure the first property we assume that
the Riemannian metric is K\"ahler with respect to some complex structure.
The second property is somewhat problematic at present, but a necessary
condition is that the K\"ahler form of the metric be in the same
de~Rham cohomology class as the K\"ahler form of some Ricci-flat metric,
and that the volume of the metric be sufficiently large.

Let ${\cal J}$ be a complex structure on $X$ for which
 the metric $g_{ij}$ is K\"ahler.
If we
pick a complex structure on $\Sigma$ which makes
its Riemannian metric K\"ahler, and which is compatible with its
orientation,  then the first term in the action can be
rewritten using the formula:
\[
\int_\Sigma \|d\phi\|^2\,d\mu=\int_\Sigma \|\bar\partial\phi\|^2\,d\mu
+\int_\Sigma \phi^*(\omega),
\]
where
$\bar\partial\phi\in\Hom(T_\Sigma^{(1,0)},\phi^*(T_{X_{\cal J}}^{(0,1)}))$
is determined by the complex structures, and where
 $\omega$ is the K\"ahler form of the metric $(g_{ij})$ on $X$.
It follows that when the classical action is evaluated on a
holomorphic map $\phi$ (i.e.,
one with $\bar\partial\phi\equiv0$), the result is simply
$$\int_\Sigma \phi^*(B+i\omega)\in\C/\Z.$$
The ``topological'' correlation functions (of both A-model and B-model
type)---when evaluated using $\sigma$-model perturbation theory---%
depend only on these extrema of the classical action, and so ultimately
 will depend only on the
choice of complex structure ${\cal J}$ and
{\em complexified K\"ahler form}\/
$\beta:=B+i\omega\in H^2(X,\C/\Z)$.

\subsection{The A-model parameter space}

Not every element of $H^2(X,\C/\Z)$ corresponds to a complexified
K\"ahler form; the ones which do, for a fixed complex structure
${\cal J}$ on $X$, constitute the
{\em complexified K\"ahler cone}
\[ {\cal K}_{\C}:=
 \{\beta\in H^2(X,\C/\Z)\,|\, \Im(\beta)
\text{ lies
within the K\"ahler cone of }X_{\cal J}\}.\]
The perturbative analysis of the $\sigma$-model is expected to be
valid in some open subset of
${\cal K}_{\C}$
containing all metrics of sufficiently large
volume, that is, in a set of the form
\[ ({\cal K}_{\C})^\circ:=
 \{\beta\in H^2(X,\C/\Z)\,|\, \Im(\beta)
\text{ lies
{\it deep}\/
within the K\"ahler cone of }X_{\cal J}\}.\]
The actual parameter space for
$\sigma$-models with complex structure ${\cal J}$ can then be represented
as $({\cal K}_{\C})^\circ/\Aut(X_{\cal J})$,
where $\Aut(X_{\cal J})$ is the
group\footnote{Typically, the group $\Aut(X_{\cal J})$
acts discretely on $({\cal K}_{\C})^\circ$.}
of diffeomorphisms of $X$ which preserve the complex structure ${\cal J}$.
We refer to this as the {\it A-model parameter space}, since the
correlation functions of the A-model are independent of the complex
structure but do depend on the parameters being described here.

A more global analysis \cite{catp,phases} reveals that the parameter space
$({\cal K}_{\C})^\circ/\Aut(X_{\cal J})$ must often be enlarged if
we wish to describe the full
moduli space of $N{=}(2,2)$
conformal field theories.  But for our present purposes, we are more
concerned with the ``large radius limit'' which occurs at the boundary of
$({\cal K}_{\C})^\circ/\Aut(X_{\cal J})$, and we need not worry about
such enlargements.

\subsection{Flat coordinates and the large radius limit}

In order to put specific coordinates on the A-model parameter
space $({\cal K}_{\C})^\circ/\Aut(X_{\cal J})$, we
need to choose a presentation for $H_2(X,\Z)$ with generators
$e_1, \dots, e_r, f_1, \dots, f_s$ and relations $m_kf_k=0$, $k=1,\dots, s$,
for some natural numbers $m_k>1$.  Thus, $f_1, \dots, f_s$ generate the torsion
subgroup,
and $e_1, \dots, e_r$ form a basis for the free abelian group
$H_2(X,\Z)/(\text{torsion})$.  We introduce the dual basis
$e^1, \dots, e^r$ of $H^2_{\text{DR}}(X,\Z)$, which will generate
the integral lattice
that provides the ambiguity in the flat coordinates.
We make the crucial assumption that
{\em each $e^j$ lies in the closure of the K\"ahler cone}.

Since $H^2(X,\C/\Z)$ is isomorphic to $\Hom(H_2(X,\Z),\C^*)$, each point
$\beta\in({\cal K}_{\C})^\circ$ can be regarded as a homomorphism, and as such
 is determined by its values on a basis, i.e., by
$q_j:=\beta(e_j)$ and
$\tau_k:=\beta(f_k)$, which must be nonzero complex numbers.
The latter are subject to the relations
$\tau_k^{(m_k)}=1$; the choice of {\em which}\/ roots of unity to use
for the $\tau_k$'s determines {\em which}\/ connected component of the
parameter space we are working with.  The $q_j$'s are exponentials of
the components of the original $2$-form, that is, when
$B+i\omega\in H^2(X,\C)/H^2_{\text{DR}}(X,\Z)$ we can write
\[B+i\omega=\frac1{2\pi i}\sum_j(\log q_j)e^j.\]
It is the logarithms $t_j:=\frac1{2\pi i}\log q_j$ which are the
``flat'' coordinates.  These are multiple-valued, and can be shifted
by independent integers.
(This
indicates how the lattice $H^2_{\text{DR}}(X,\Z)$ specifies the ambiguity
in the flat coordinates.)
 However, the corresponding vector fields and $1$-forms
\[\frac{\partial}{\partial t_j}=2\pi i\, q_j\,\frac{\partial}{\partial q_j}
\qquad \text{and} \qquad
dt_j = \frac1{2\pi i}\,\dlog q_j
\]
are single-valued.

In order to study the large radius limit, we restrict our attention to
those K\"ahler classes which lie in the cone
\[{\cal C}={\cal C}_{\vec{e}}:=\{\omega=\sum\omega_je^j\,|\, \omega_j>0\}\]
spanned by the chosen basis vectors.
If the action of $\Aut(X_{\cal J})$ on ${\cal K}_{\C}$ is discrete,
then it will be possible to find such bases with the property that
${\cal C}_{\vec{e}}$ is disjoint from its translates under
$\Aut(X_{\cal J})$.  (In any case, we shall ignore the action
of $\Aut(X_{\cal J})$
for the time being.)
The corresponding complexified cone
\[{\cal C}_{\C}:=\{\beta\in H^2(X,\C/\Z)\,|\, \Im(\beta)\in{\cal C}\}
\subset {\cal K}_{\C}\]
is described in coordinates by the condition
\[\Im(\frac1{2\pi i}\log q_j)=-\frac1{2\pi}\log|q_j|>0\text{ for all }j,\]
or equivalently,
\[0<|q_j|<1\text{ for all }j.\]
To find the large radius limit, we should rescale $\omega\to \lambda\omega$,
and let $\lambda$ grow to infinity.  Under such a rescaling, we have
\[q_j\mapsto |q_j|^{(\lambda-1)}\,q_j .\]
Thus, all points in ${\cal C}_{\C}$ flow towards $q_j=0\ \forall j$ under
this
rescaling, and $q_j=0\ \forall j$ should be taken as the ``large radius
limit.''\footnote{It is not yet precisely clear how one should interpret
the torsion variables $\tau_k$ in the large radius limit,
but see \cite{stable} for some steps in this direction.
For the purposes of this paper, we work with the component in which
$\tau_k=1$ for all $k$.}
We form a partial compactification of our parameter space by enlarging it
to include all $q$'s such that
\[0\le |q_j|<1\text{ for all }j.\]
On this enlarged space, the $q_j$'s occur as natural coordinates, and
the ``boundary'' of the space is a divisor with normal crossings.

There is a natural identification which can be made between the space
of marginal operators for the A-model and the vector fields
$\partial/\partial t_j$ corresponding to the flat coordinates.
When we calculate three-point functions with respect to these coordinates,
we find an expansion of the form
\[\big\langle
\frac{\partial}{\partial t_j}
\frac{\partial}{\partial t_k}
\frac{\partial}{\partial t_\ell}
\big\rangle =
e^j \cup e^k \cup e^\ell |_{[X]} + O(q),
\]
where $O(q)$ represents the instanton corrections to the classical value,
which contain the data about the enumeration of rational curves.
If we write this in terms of the (single-valued)
 coordinates $q_j$ at the large radius limit point,
we find
\[(2\pi i)^3\,q_jq_kq_\ell\,
\big\langle
\frac{\partial}{\partial q_j}
\frac{\partial}{\partial q_k}
\frac{\partial}{\partial q_\ell}
\big\rangle =
e^j \cup e^k \cup e^\ell |_{[X]} + O(q).
\]
In other words, the three-point function has poles along the boundary
divisor in the $q$-coordinates.  Moreover, the leading order term in
a Laurent expansion of a three-point function
 picks out the corresponding cohomological quantity.

The analysis we have given depends on a choice of basis; we defer to
section \ref{s:ambiguity}
a discussion of what happens when the basis is changed.

\section{Maximally unipotent monodromy} \label{s:maxunip}

The structure which we have found in the A-model---a partial compactification
of the parameter space, with poles of the correlation functions along the
compactification divisor---will now serve as a guide to making enumerative
predictions by means of B-model calculations.  In order to carry
this out, we must the analyze compactifications of the B-model moduli space.

Given an arbitrary compactification of the B-model moduli space, we are
always free to blow   up along the boundary until the boundary becomes
 a normal crossings divisor.  The only remaining singularities of the space
after such a blowup would lie in the interior of the moduli space (where there
may well be quotient singularities associated with complex structures for
which the automorphism group is larger than generic).
Even those can be removed by passing to a finite cover.

When the boundary is a normal crossings divisor,  the monodromy theorem
\cite{monodromy} guarantees that the monodromy of the periods around
each component of the boundary is a quasi-unipotent transformation
(unipotent after passing to a finite cover).
Unipotent monodromy appears
to be necessary in order to correctly reproduce the behavior of the A-model.
We will therefore assume that the monodromy transformations near the
points we seek
are in fact unipotent.

In order to analyze the boundary in detail and search for the mirrors of large
radius limit points, we  restrict our attention
to a local situation in which a product of
punctured disks (with coordinates $s_j$) has been embedded in the interior
of our moduli space in such a way that the limit points $s_j\to0$ are
mapped to the boundary.
Let $T^{(j)}$ be the monodromy transformation about the $j^{\text{th}}$
coordinate (counterclockwise), with respect to some fixed basepoint $P$
near the origin.

For discussions of monodromy, it is more convenient to
 represent each period $\int_{\gamma}\Omega(s)$ by means of cup product with
a cohomology class $g\in H^3(Y_P,\C)$, i.e.,
\[\int_{\gamma}\Omega(s)=\langle
g\,|\,\Omega(s)\rangle:=\int_{Y_P}g\wedge\Omega(s).\]
The cycle $g$ extends to a multi-valued
section of the local system $R^3\pi_*\C_{\cal Y}$,
and the corresponding period
is also multi-valued.  However, according to the nilpotent orbit
theorem \cite{schmid}, the section
\[\exp(-\frac1{2\pi i}\log s\cdot\log T)\,g\in\Gamma({\cal V})\]
of the locally free sheaf ${\cal V}=R^3\pi_*\C_{\cal Y}
\otimes{\cal O}_S$
is single-valued.  We introduce\footnote{This sign convention differs
from \cite{compact}, but agrees with \cite{Deligne}.}
 $N^{(j)}=-\log T^{(j)}$ so that the corresponding single-valued section
can be written as
\[\widetilde{g}:=\exp(\frac1{2\pi i}\sum(\log s_j)N^{(j)})\,g.\]

We now consider the conditions on periods needed to match the
behavior of the A-model.  First, the period $\int_{\gamma_0}\Omega(s)
=\langle g^0\,|\,\Omega(s)\rangle$ should be
single-valued, so we need to find a cycle $g^0$ such that $N^{(j)} g^0=0$
for all $j$.  Second, the monodromy on the period
 $\int_{\gamma_j}\Omega(s)
=\langle g^j\,|\,\Omega(s)\rangle$
should only involve $\int_{\gamma_0}\Omega(s)$, in order that the
multi-valuedness of the ratios shifts them by constants.  Thus, we
need for $N^{(j)} g^k$ to be a multiple of $g^0$ for every $j$ and $k$.
In fact, if we write
\[ N^{(j)}g^k=m^{jk}g^0,\]
then the matrix $(m^{jk})$ must be invertible in order to solve for
 coordinates with the desired monodromy properties.
Letting $(m_{k\ell})$ denote the inverse matrix of $(m^{jk})$,
we have
\[(T^{(j)}-I)(-\sum g^\ell m_{\ell k})=\sum N^{(j)}g^\ell m_{\ell
k}=\delta^j_k\,g^0,\]
so this cohomology class ``$-\sum g^\ell m_{\ell k}$''
determines the ratio of periods which has the correct monodromy
properties.
The corresponding (multi-valued) coordinates are then
\[
\frac1{2\pi i}\,\log z_k=\frac{-1}{\langle g^0\,|\,\Omega\rangle}
\sum_{\ell=1}^r\langle g^\ell\,|\,\Omega\rangle m_{\ell k}.
\]
 As in the case of
the A-model, by exponentiating these we obtain single-valued coordinates $z_k$
which
extend across the boundary.

The coordinates as written are not uniquely specified, and in fact we have
not yet used all of the information available to us by comparison with the
A-model.  What we have not yet considered is the three-point functions of
the B-model, which should be mirror to the three-point functions of the
A-model.
In order to describe these, we must fix a particular choice $\Omega(s)$ of
 holomorphic $3$-forms on the fibers, which can be thought of as fixing
the gauge in the bundle $\pi_*\omega_{{\cal X}/S}$ whose fibers are
the spaces of global holomorphic $3$-forms on the fibers of $\pi$.
Moreover,
as in the case of the A-model, the marginal operators whose correlation
functions we wish to calculate can be identified with vector fields on
the moduli space.  Any system of  coordinates $s_j$ has an associated
collection of vector fields $\partial/\partial s_j$, and with respect
to these, the three-point function can be written
\[\big\langle
\frac{\partial}{\partial s_j}
\frac{\partial}{\partial s_k}
\frac{\partial}{\partial s_\ell}
\big\rangle :=
\int_Y\Omega(s)\wedge\nabla_{s_j}\nabla_{s_k}\nabla_{s_\ell}\Omega(s)\]
where $\nabla_{s_j}\varphi$ represents the directional derivative
$(\nabla\varphi)\lhk \frac{\partial}{\partial s_j}$.

When we calculate these three-point functions near the boundary of
the moduli space, we should expect to find poles (in order to replicate
the behavior of the A-model moduli space), and indeed the presence
of poles in the extension of $\nabla$ to the boundary is a well-known
phenomenon in algebraic geometry (cf.\ \cite{regsings}).
These poles arise from the behavior
of the single-valued sections of ${\cal V}$ under differentiation:
if we calculate using the single-valued section $\widetilde g$ introduced
above, we find
\[\nabla(\widetilde g)=\frac1{2\pi i}\sum(\dlog s_j)N^{(j)}\,\widetilde g.\]
The coefficients $\dlog s_j$ in this expression are $1$-forms with poles
along the boundary.

Taking three directional derivatives $\nabla_{s_j}\nabla_{s_k}\nabla_{s_\ell}$,
we see that the leading term
in a Laurent expansion has coefficient proportional to something of the
form
\[\left(\frac1{2\pi i}\right)^3 N^{(j)}N^{(k)}N^{(\ell)}\,\widetilde g .\]
Now the comparison with the A-model tells us that at least some of these
coefficients must be nonzero, since they should reproduce the intersection
numbers $e^j\cup e^k\cup e^\ell|_{[X]}$ of the mirror partner.
(In fact, Poincar\'e duality on $X$ implies some rather strong conditions
on these intersection numbers, which must be replicated by the
coefficients we are calculating on the B-model side.)
The simple fact that any of these numbers is nonzero, though, immediately
implies that the order of unipotency of the monodromy transformations is
in some sense ``maximal''.  (Any quartic expressions in the $N$'s must
vanish for dimension reasons, and so cubic expressions are the maximal
possible degree for a non-vanishing expression.)
Furthermore, the fact that the $N^{(j)}$'s define a limiting
mixed Hodge structure in which $h^{(3,0)}=1$ implies that the images
of all of the cubic expressions lies in a one-dimensional space $W_0$.

Then, by including
the relations deduced from Poincar\'e duality, we arrive at the following
definition.\footnote{The original definition given in \cite{mirrorguide}
only applied to the one-parameter case; this was
 extended to several parameters by Deligne \cite{Deligne}
and  the author \cite{compact}. (We follow the version given in
\cite{compact}.)}

\begin{definition}
A normal crossings boundary  point $P$ of $S$
is called a\/ {\em maximally unipotent point} under the following conditions.
\begin{enumerate}
\item
$P$ lies at the intersection of $r=\dim S$ local boundary components $B_j$,
and the
 monodromy transformations $T^{(j)}$ around these components
 are all unipotent.
\item
Let $N^{(j)}=-\log T^{(j)}$, let $N :=\sum a_jN^{(j)}$ for some $a_j>0$,
and define
\begin{eqnarray*}
W_0&:=&\Im(N ^3)\\
W_1&:=&\Im(N ^{2})\cap\Ker N \\
W_2&:=&\left(\Im(N )\cap\Ker(N )\right)+\left(
\Im(N ^{2})\cap\Ker(N ^2)\right).
\end{eqnarray*}
Then $\dim W_0=\dim W_1=1$ and $\dim W_2=1+\dim(S)$.
\item
Let $g^0,g^1,\dots,g^r$ be a basis of $W_2$ such that $g^0$ spans $W_0$,
and define $m^{jk}$ by
$N^{(j)}g^k=m^{jk}g^0$ for $1\le j,k\le r$.
Then $m:=(m^{jk})$ is an invertible matrix.
\end{enumerate}
(The spaces $W_0$ and $W_2$ are independent of the choice of coefficients
$\{a_j\}$
{\rm\cite{CK,deligne:weil2}},
and the invertibility of $m$ is independent of the choice
of basis $\{g^k\}$.)
\end{definition}

When we restrict to maximally unipotent boundary points, the single-valued
$1$-forms
\[
\frac1{2\pi i}\,\dlog z_k=d\left(\frac{-1}{\langle g^0\,|\,\Omega\rangle}
\sum_{\ell=1}^r\langle g^\ell\,|\,\Omega\rangle m_{\ell k}\right).
\]
are independent of the choice of basis $\{g^k\}$ and of $3$-form $\Omega$.
The coordinates $z_k$ themselves {\em do}\/ depend on the choice of basis,
but only through multiplicative constants:  a change of basis replacing
$g^k$ by $\sum_{\ell=0}^kc^k_\ell g^\ell$ will induce
\[\frac1{2\pi i}\log z_k \mapsto \frac{c^k_0}{c^0_0} + \frac1{2\pi i}\log z_k\]
and so
\[z_k\mapsto e^{2\pi i(c^k_0/c^0_0)}\,z_k.\]
Determining the ``constants of integration'' which specify $z_k$ once $\dlog
z_k$
is known is the most delicate part of finding the mirror map.  We will address
this issue in section \ref{s:mirrormap}.

In practice,
 computing the monodromies about all boundary components and locating
which points on the boundary have maximally unipotent monodromy is a
challenging task.  In fact, this computation has only been carried out
fully in a few examples \cite{CdGP,2param1,2param2}.  In the special cases
of complete intersections in toric varieties, there is another method
which has been use to locate such points: one finds the natural ``toric''
limit points in the toric moduli space which correspond to K\"ahler cones
of possible birational models of the mirror, and these turn out to have
maximally unipotent monodromy (as follows from \cite{Bat:qcoho,BK}).

This alternate method must be used with some caution, for {\em it is possible
to have maximally unipotent boundary points which are not toric boundary
points.}  An explicit example of this phenomon was seen in \cite{2param2},
where there is a non-toric boundary point with maximally unipotent monodromy.
Interestingly, in that example there is also an additional (non-toric)
 discrete symmetry
of the toric moduli space by which one must quotient  in order to obtain the
true B-model moduli space.  That additional discrete symmetry
exchanges the toric and non-toric points with maximally unipotent monodromy.
It would be interesting to know whether or not this is true in general:
in the toric complete intersection case, given a boundary point with maximally
unipotent monodromy, does there always exist a discrete symmetry of the
moduli space which maps this point to a toric boundary point?

\section{Equivalence among boundary points} \label{s:ambiguity}

Our discussion in section \ref{s:largeradius} of coordinates near the
large radius limit depended on the choice of basis for
$H_2(X,\Z)/(\text{torsion})$,
or equivalently, on the choice of simplicial cone ${\cal C}\subset{\cal K}$.
The fact that different choices of simplicial cone
(always contained in the K\"ahler cone) lead to apparently different ``large
radius
limit'' points in the A-model parameter space
should not be too surprising.  The limit point we seek
is actually a boundary point of our parameter space, and what
we are finding is that there are different ways to compactify the
space.  Since we are using compactifications which (locally) have the
structure of an algebraic variety, we should expect birational
modifications along the boundary to provide a mechanism for passing
between compactifications and indeed that is what happens with our
choice of cones.  Subdividing a given cone into smaller ones precisely
corresponds to blowing up, as in toric geometry.

For example, if we start from a basis $e^1, \dots, e^r$ and blowup the origin
in
the coordinate chart $(q^1,\dots,q^r)$, we
find  new coordinate charts
after the blowup, with coordinates
$(q_1,\frac{q_2}{q_1},\dots,\frac{q_r}{q_1})$,
\dots, $(\frac{q_1}{q_r},\dots,\frac{q_{r-1}}{q_r},q_r)$, respectively.
(The corresponding bases are $\{e^1+\cdots+e^r, e^2,\dots,e^r\}$,
\dots,  $\{e^1,\dots,e^{r-1}, e^1+\cdots+e^r\}$, respectively.)
Rescaling the metric and taking $\lambda\to\infty$ sends $(q_1,\dots,q_r)$
to the origin in the first chart when $1>|q_1|>|q_j|$ ($\forall\ j\ne1$),
sends it to the origin in the second chart  when $1>|q_2|>|q_j|$ ($\forall\
j\ne2$),
and so on. All of these ``origins'' can thus
lay claim
to being ``the large radius limit'' associated to at least {\it part}\/
of the A-model parameter space.

Conversely, if we have a partial
compactification of the A-model parameter space
which includes more than one large radius limit point (each associated
with a different basis $e^1,\dots,e^r$, and with a different domain inside
the moduli space), we should attempt to blow down this  space
 to produce a partial compactification with a
 single large radius limit point for the entire moduli space.
These blowdowns are similar to those arising in toric geometry,
and will often lead to singularities in the compactified space.
The instanton contributions to correlation functions are still suppressed
in such a limit, in spite of the singularities---we must accept the
possibility that the ``true'' large radius limit point is not a smooth
point.

(Note that all of the large radius limit points under discussion are
associated to a  single K\"ahler cone.  It is also possible
to consider other large radius limit points associated to the
K\"ahler cones of  different birational models of $X$.
This leads to topology-changing transitions \cite{catp},
and we would not expect to collapse  those limit points
to a single point by blowing down.)

Comparison between different cones can be accomplished by considering
the canonical $1$-forms $\dlog q_j$.  These are intrinsically defined,
and should only change by a constant change-of-basis matrix when
moving from one large radius limit point to another (within the same
K\"ahler cone).  These $1$-forms will therefore define a local
system $\cal L$ in a neighborhood of all of the exceptional divisors
of a potential blowing-down map associated to the K\"ahler cone.

We are thus led to introduce an equivalence relation among boundary points
of the A-model parameter space:
two boundary points $P$ and $Q$ are equivalent if there exists a connected
subset $\Xi$ of the boundary  containing both points and a local system
${\cal L}$ defined in a neighborhood of $\Xi$ which is spanned by
the canonical $1$-forms $\dlog q_j$ at any maximally unipotent point within
$\Xi$.  (For further details about this construction, we refer the
reader to \cite{compact}.)

Even when we expect to be able to
 blow down and are willing to allow singularities,
it may prove to be impossible to perform the desired blowing down,
due to the presence
of an infinite number of large radius limit points.
We describe this phenomenon in an explicit example,
following  \cite{where}.  Suppose that  the K\"ahler cone is described as
$\frac2{1-\sqrt5}\,y<x<\frac2{1+\sqrt5}\,y$. Then (as shown in figure 1)
attempting to cover the cone using integral bases leads to a
sequence of rays with slopes
$\frac y x = \dots,-\frac58,-\frac23,-1,\frac10,2,\frac53,\frac{13}8,\dots$
which asymptotically approach the
walls\footnote{The figure does not include these walls---the
 limiting rays with
irrational slope $\frac{1\pm\sqrt5}2$ ---since they are less
than a line-width's distance from the outer rays as shown (at the
level of resolution of the figure).} of the cone.
Each adjacent pair of rays in the sequence
gives rise to a distinct large radius limit
point.

$$\vbox{\xpt
\baselineskip=12pt
\iffigs
\centerline{\epsfxsize=8cm\epsfbox{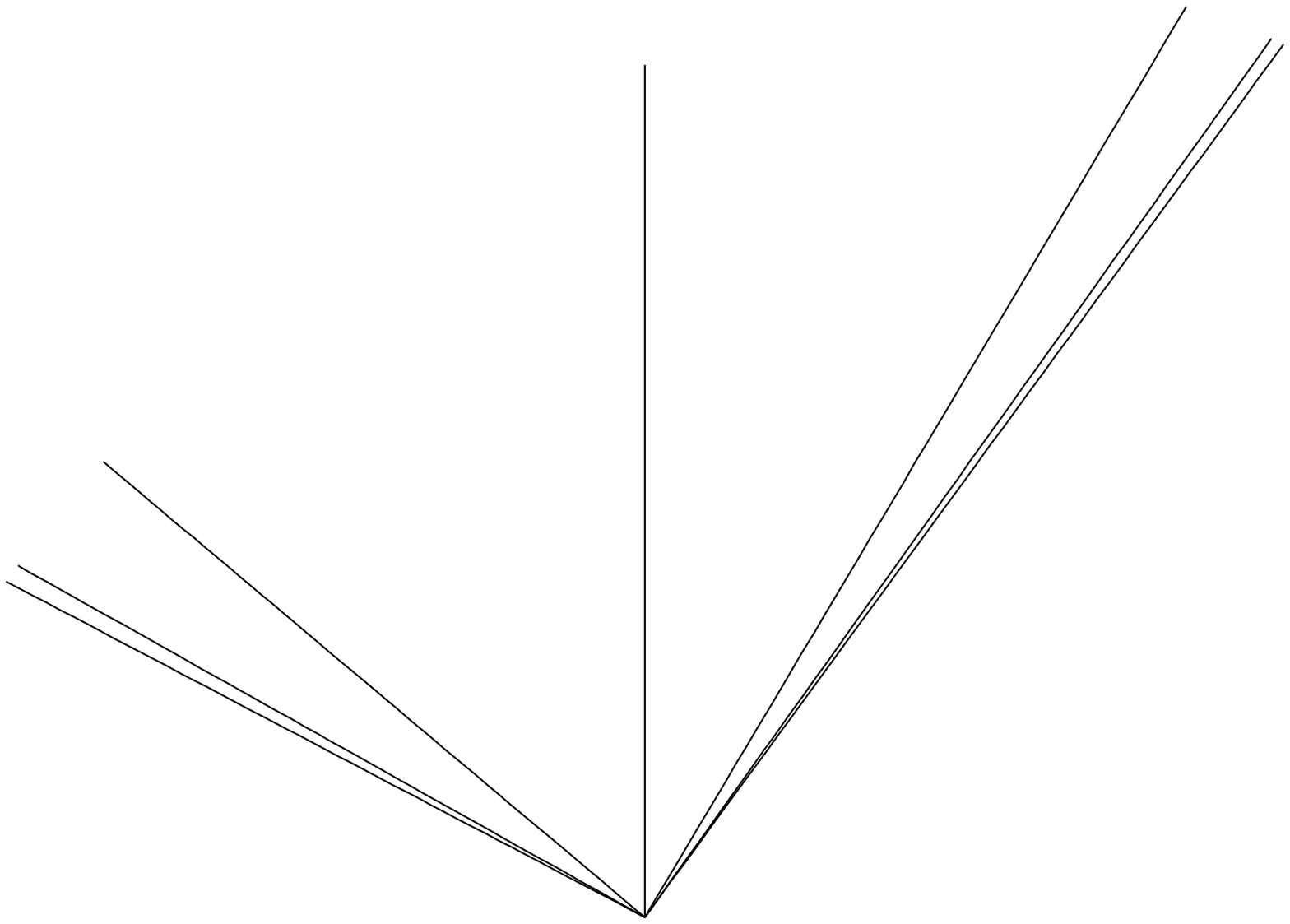}}
\else \vglue1cm \fi
\centerline{Figure 1. Decomposing the cone
$\frac2{1-\sqrt5}\,y<x<\frac2{1+\sqrt5}\,y$.}}$$

However, when we include the action of $\Aut(X_{\cal J})$ in our analysis,
it may become
possible to do the blowdowns---an infinite number of large radius limit
points may turn into a finite number
after these identifications \cite{compact,kcone}.
In the example above, an automorphism acting on the cone as
$(x,y)\mapsto(2x+3y,3x+5y)$ leads from an infinite number of large radius
limit points on the original K\"ahler moduli space to two remaining
points on the quotient space.
There are two boundary divisors (after taking the quotient), and they
meet in two large radius limit points.
The quotient space can then be blown down
explicitly using methods of Hirzebruch \cite{hirzebruch},
leading to a  surface singularity with local
equation $w^2=(u^3-v^2)(u^2-v^3)$.  This is illustrated in figure 2.
In general, $\Xi$ will be a subset of the compactified moduli space
only after taking such a quotient, which is why we use a local system ${\cal
L}$
rather than simply a trivial sheaf.

$$\vbox{\xpt
\baselineskip=12pt
\iffigs
\centerline{\epsfxsize=10cm\epsfbox{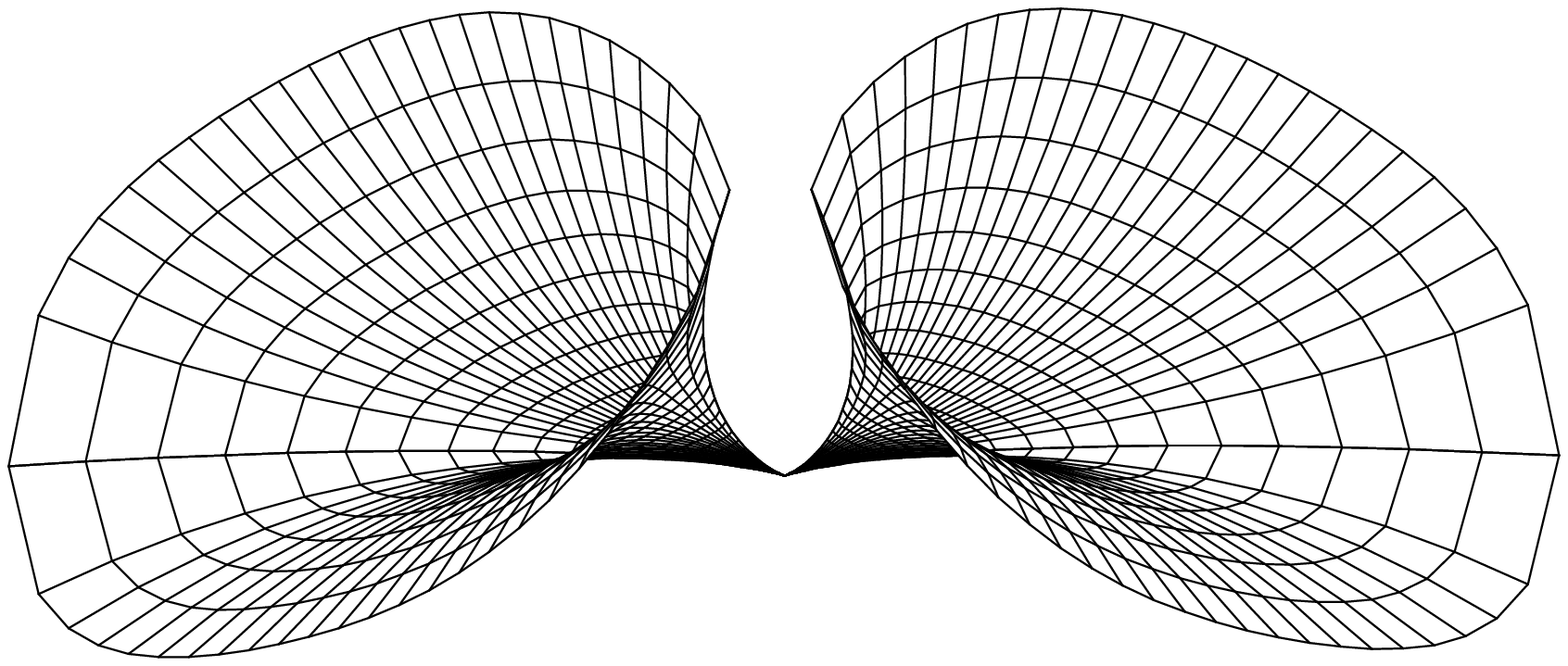}}
\else\vglue1cm\fi
\centerline{Figure 2. The blown down moduli space $w^2=(u^3-v^2)(u^2-v^3)$}}
$$

Applying these ideas to the analysis of the B-model moduli space,
we see that we should consider two large radius limit points to be
equivalent
 when there exists a connected
subset $\Xi$ of the boundary  containing both points and a local system
${\cal L}$ defined in a neighborhood of $\Xi$ which is spanned by
the canonical $1$-forms $\dlog z_k$ at any maximally unipotent point within
$\Xi$.  This can be effectively computed if we know the mirror map at
each maximally unipotent point; in fact, it is enough to calculate the
leading terms in Laurent expansions of $3$-point functions.  Further
details are in \cite{compact}.

\section{Determining the mirror map (two conjectures)} \label{s:mirrormap}

As pointed out in
section \ref{s:maxunip}, the most delicate part of determining
the mirror map is specifying the constants of integration, passing from
the canonical $1$-forms $\dlog z_j$ to the actual multi-valued coordinates
$z_j$ themselves.  There are two conjectures which have been used to
determine these constants.  The first one (stated here in detail for the first
time) is completely general, applying in principle to all Calabi--Yau
threefolds, while the second one is special to the case of toric hypersurfaces.

In addition to specifying the coordinates, one would like to specify a
particular choice of holomorphic $3$-forms $\Omega(z)$ in order to determine
the $3$-point functions precisely.  This is something which
the conjectures will also do.

\subsection{A conjecture about integral cohomology}

The first conjecture for determining the mirror map involves the integral
cohomology groups of a Calabi--Yau threefold.  The conjecture is quite
natural,  and there is a bit of
evidence for it in a few specific examples.  We will give additional evidence
here.

Simply put, we conjecture that the canonical coordinates and canonical
gauge for $\Omega(z)$ should be given by periods over {\em integer-valued}\/
cohomology classes, in the following precise sense.
The classes $g^0$, $g^1$, \dots, $g^r$
should be chosen from $H^3(Y,\Z)$, in such a way that
$g^0$ spans $W_0\cap H^3(Y,\Z)$ and the entire set $g^0$, $g^1$, \dots, $g^r$
spans $W_2\cap H^3(Y,\Z)$.  We conjecture that if this is done, and if
we write $N^{(j)}g^k=m^{jk}G^{(0)}$ as in the previous section, then
the matrix $(m^{jk})$ {\em is invertible over the integers}.  If this
is true, then the mirror map will be uniquely specified by using the
integral periods, and  the gauge $\Omega(z)$
for which $\langle g^0\,|\,\Omega(z)\rangle=\pm1$ will be uniquely
specified (up to sign) as well.
(This is because in any change of basis
$g^k\mapsto\sum_{\ell=0}^kc^k_\ell g^\ell$ preserving the integral
structure, we will have $c^0_0=\pm1$ and $c^k_0\in\Z$ so that
$\exp(2\pi i(c^k_0/c^0_0))=1$.)

Our conjecture is motivated by the observation that the integral structure
on $H^2(X)$ controls the choice of coordinates there.
In fact, the action of $N^{(j)}$ on $H^3(Y)$ can be seen as mirroring
the action ``cup product with $e^j$'' on $H^*(X)$, where $e^j\in H^2(X)$
is an integral class, part of the basis determining the coordinates.
The idea that the integral structure on $H^3$ should mirror the integral
structure on $H^0\oplus H^2\oplus H^4\oplus H^6$ is not a new one---it
was explicitly mentioned in \cite{AL:qag}, for example, and it was
used implicitly in the calculations of \cite{CdGP}
(cf.\ also \cite{Cd2}).  There is not
a lot of evidence for this equivalence, however, other than the examples which
we
describe here.

As a practical matter, our conjecture can be tested in the following
way.  Compute the monodromy matrices $T^{(j)}$ with respect to a basis
of integral cohomology.  There must then be a rank one matrix $M$
(with image $W_0$) such that
$N^{(j)}N^{(k)}N^{(\ell)}=c^{jk\ell}M$ for all $j$, $k$, $\ell$, where
$c^{jk\ell}:=e^j\cup e^k\cup e^\ell|_{[X]}$ are the intersection
numbers on the mirror partner $Y$ of $X$.
The conjecture states that $M$ should
be a primitive
integral matrix.\footnote{In this version of testing the conjecture, it is
assumed
that a mirror partner is known.  One could test the conjecture without this
assumption by finding the primitive integral matrix $M$ first, calculating
the corresponding coefficients $c^{jk\ell}$, and checking to see if
they have the numerical properties compatible with Poincar\'e duality
over the integers.}

The fact that the mirror map can perhaps be completely
determined by looking at the integral
structure was pointed out in \cite{mirrorguide}, where it
was shown that this conjecture holds for the case of the quintic-mirror
(using calculations from \cite{CdGP}), and that the integral basis leads
to the correct mirror map.
In the two-parameter examples of \cite{2param1} the same principle was
used to determine the mirror map, with equal success.\footnote{The rank one
matrices in those papers---denoted by
$Y$ in eq.~(7.4) of \cite{2param1} and also by $Y$ preceding eq.~(6.4) in
\cite{2param2}---have the property of being primitive integral matrices,
although this was not pointed out in those papers.}  We will give
additional evidence for our conjecture by verifying it (and checking
that it produces the ``correct'' mirror map)
 in the three further one-parameter
examples studied in \cite{picard-fuchs,font,kt1}.

We must repeat the verification made in Appendix C of \cite{mirrorguide}
that the monodromy of the actual period functions has a certain form.
In fact, we will find a somewhat better normalization of that form
this time.  We will use the explicit  monodromy calculations
from \cite{kt1};  a similar calculation could be done
using \cite{font} (which uses a different normalization of the parameter,
making it difficult to compare to the present approach).

Our verification is displayed in table \ref{tab:mat}.
We show in the second column of the table the monodromy
matrix $A$ as calculated in \cite{kt1}.  That matrix was not
quite uniquely specified by the data with which those authors were working.
In particular, there is freedom to replace $A$ by $A'=m'A(m')^{-1}$ for
any matrix $m'\in\text{Sp}(4,\Bbb Z)$ whose second and fourth rows are
the same as that of the identity matrix.  We make a choice of $m'$,
shown in the third column of the table, and calculate $A'$ in the
fourth column.  Notice that the result takes the form
\[A'=\begin{bmatrix}
\phantom{-}1 & -1             & \phantom{-}0 & \phantom{-}1 \\
\phantom{-}0 & \phantom{-}1   & \phantom{-}0 & -1 \\
-\lambda     & \phantom{-}0   & \phantom{-}1 & \phantom{-}0 \\
-\lambda     & \phantom{-}\mu & \phantom{-}1 & \phantom{-}{1-\mu}
\end{bmatrix}\]
where $(\lambda,\mu)$ are as given in the fifth column of the table.

With $A'$ in the given form, we can calculate the monodromy around infinity
as
\[T_\infty:=T^{-1}\,{A'}^{-1}=\begin{bmatrix}
\phantom{-}1 & \phantom{-}1 & \phantom{-}0 & \phantom{-}0 \\
\phantom{-}0 & \phantom{-}1 & \phantom{-}0 & \phantom{-}0 \\
\phantom{-}\lambda & \phantom{-}\lambda & \phantom{-}1 & \phantom{-}0 \\
\phantom{-}0 & -\mu & -1 & \phantom{-}1
\end{bmatrix},\]
where
\[T:=\begin{bmatrix}1&0&0&0\\0&1&0&1\\0&0&1&0\\0&0&0&1\end{bmatrix}.\]
We can thus easily see that $T_\infty$ is unipotent, with $(T_\infty-I)^4=0$.

It is then a straightforward computation to see that
\[(-\log T_\infty)^3
=\begin{bmatrix}
\phantom{-}0 & \phantom{-}0 & \phantom{-}0 & \phantom{-}0 \\
\phantom{-}0 & \phantom{-}0 & \phantom{-}0 & \phantom{-}0 \\
\phantom{-}0 & \phantom{-}0 & \phantom{-}0 & \phantom{-}0 \\
\phantom{-}0 & \phantom{-}\lambda & \phantom{-}0 & \phantom{-}0
\end{bmatrix}.\]
Note that $\lambda$ is precisely the triple-self-intersection of an
integral generator of $H^2$ of the mirror partner, verifying the
conjecture in these cases.

\begin{table}
\begin{center}
$\begin{array}{|c|c|c|c|c|} \hline
\vphantom{\Big(}
k & A & m' & A'=m'A(m')^{-1} & (\lambda,\mu) \\
\hline
&&&&\\
5 &
\begin{bmatrix}
-9 & -3 & \phantom{-}5 & \phantom{-}3 \\
\phantom{-}0 & \phantom{-}1 & \phantom{-}0 & -1 \\
-20 & -5 & \phantom{-}11 & \phantom{-}5 \\
-15 & \phantom{-}5 & \phantom{-}8 & -4
\end{bmatrix}
&
\begin{bmatrix}
\phantom{-}2 & \phantom{-}0 & -1 & \phantom{-}0 \\
\phantom{-}0 & \phantom{-}1 & \phantom{-}0 & \phantom{-}0 \\
-5 & \phantom{-}0 & \phantom{-}3 & \phantom{-}0 \\
\phantom{-}0 & \phantom{-}0 & \phantom{-}0 & \phantom{-}1
\end{bmatrix}
&
\begin{bmatrix}
\phantom{-}1 & -1 & \phantom{-}0 & \phantom{-}1 \\
\phantom{-}0 & \phantom{-}1 & \phantom{-}0 & -1 \\
-5 & \phantom{-}0 & \phantom{-}1 & \phantom{-}0 \\
-5 & \phantom{-}5 & \phantom{-}1 & -4
\end{bmatrix}
& (5,5) \\
&&&&\\
\hline
&&&&\\
6 &
\begin{bmatrix}
\phantom{-}1 & -1 & \phantom{-}0 & \phantom{-}1 \\
\phantom{-}0 & \phantom{-}1 & \phantom{-}0 & -1 \\
-3 & -3 & \phantom{-}1 & \phantom{-}3 \\
-6 & \phantom{-}4 & \phantom{-}1 & -3
\end{bmatrix}
&
\begin{bmatrix}
\phantom{-}1 & \phantom{-}0 & \phantom{-}0 & \phantom{-}0 \\
\phantom{-}0 & \phantom{-}1 & \phantom{-}0 & \phantom{-}0 \\
-3 & \phantom{-}0 & \phantom{-}1 & \phantom{-}0 \\
\phantom{-}0 & \phantom{-}0 & \phantom{-}0 & \phantom{-}1
\end{bmatrix}
&
\begin{bmatrix}
\phantom{-}1 & -1 & \phantom{-}0 & \phantom{-}1 \\
\phantom{-}0 & \phantom{-}1 & \phantom{-}0 & -1 \\
-3 & \phantom{-}0 & \phantom{-}1 & \phantom{-}0 \\
-3 & \phantom{-}4 & \phantom{-}1 & -3
\end{bmatrix}
& (3,4) \\
&&&&\\
\hline
&&&&\\
8 &
\begin{bmatrix}
\phantom{-}1 & -1 & \phantom{-}0 & \phantom{-}1 \\
\phantom{-}0 & \phantom{-}1 & \phantom{-}0 & -1 \\
-2 & -2 & \phantom{-}1 & \phantom{-}2 \\
-4 & \phantom{-}4 & \phantom{-}1 & -3
\end{bmatrix}
&
\begin{bmatrix}
\phantom{-}1 & \phantom{-}0 & \phantom{-}0 & \phantom{-}0 \\
\phantom{-}0 & \phantom{-}1 & \phantom{-}0 & \phantom{-}0 \\
-2 & \phantom{-}0 & \phantom{-}1 & \phantom{-}0 \\
\phantom{-}0 & \phantom{-}0 & \phantom{-}0 & \phantom{-}1
\end{bmatrix}
&
\begin{bmatrix}
\phantom{-}1 & -1 & \phantom{-}0 & \phantom{-}1 \\
\phantom{-}0 & \phantom{-}1 & \phantom{-}0 & -1 \\
-2 & \phantom{-}0 & \phantom{-}1 & \phantom{-}0 \\
-2 & \phantom{-}4 & \phantom{-}1 & -3
\end{bmatrix}
& (2,4) \\
&&&&\\
\hline
&&&&\\
10 &
\begin{bmatrix}
\phantom{-}1 & \phantom{-}0 & \phantom{-}1 & \phantom{-}0 \\
\phantom{-}0 & \phantom{-}1 & \phantom{-}0 & -1 \\
\phantom{-}0 & \phantom{-}1 & \phantom{-}1 & -1 \\
\phantom{-}1 & \phantom{-}3 & \phantom{-}1 & -2
\end{bmatrix}
&
\begin{bmatrix}
\phantom{-}0 & \phantom{-}0 & -1 & \phantom{-}0 \\
\phantom{-}0 & \phantom{-}1 & \phantom{-}0 & \phantom{-}0 \\
\phantom{-}1 & \phantom{-}0 & \phantom{-}0 & \phantom{-}0 \\
\phantom{-}0 & \phantom{-}0 & \phantom{-}0 & \phantom{-}1
\end{bmatrix}
&
\begin{bmatrix}
\phantom{-}1 & -1 & \phantom{-}0 & \phantom{-}1 \\
\phantom{-}0 & \phantom{-}1 & \phantom{-}0 & -1 \\
-1 & \phantom{-}0 & \phantom{-}1 & \phantom{-}0 \\
-1 & \phantom{-}3 & \phantom{-}1 & -2
\end{bmatrix}
& (1,3) \\
&&&&\\
\hline
\end{array}$
\end{center}

\medskip

\caption{Monodromy calculations}
\label{tab:mat}
\end{table}

\subsection{The monomial-divisor mirror map}

In the case of toric hypersurfaces, there is an alternate conjectural method
for
specifying the mirror map, proposed in \cite{mondiv},\footnote{Some
signs were left unspecified in \cite{mondiv}. The proposal for determining
the signs  given in \cite{small} is now in doubt; an alternate proposal
\cite{summing} has much evidence in its favor.}
and used with great success in \cite{catp,small,HKTY1}
(see also \cite{2param1,2param2}).  Briefly, the parameters
on both the A-model and B-model sides can be described by remarkably
similar combinatorics; this similarity is used to write a conjecture
for the derivative of the mirror map, which specifies the constants of
integration.  The conjecture was extended in \cite{summing} to also
specify the ``algebraic gauge'' which should be used as a starting point
for determining the natural gauge $\Omega(z)$; in addition, much evidence
was amassed in \cite{summing} in favor of this approach.  We refer the
reader to \cite{mondiv} and \cite{summing} for details.

\section{Making enumerative predictions} \label{s:predictions}

We are finally ready to put together all of the ingredients and describe
the process of making enumerative predictions.  The things which we are
going to predict are the ``numbers of rational curves'' on a Calabi--Yau
threefold, in the precise form of the ``Gromov--Witten invariants'' of
the threefold.  A mathematical version of these invariants has been extensively
investigated \cite{MS,RT1}  using Gromov's symplectic geometry techniques
 \cite{Gromov}
which had inspired Witten's original work on the invariants
\cite{tsm,W:aspects}.
(An alternate proposed definition purely within algebraic geometry is
currently under development by Kontsevich and Manin \cite{KM,Kontsevich}.)

The steps in an enumerative prediction are these: given a proposed mirror pair
$(X,Y)$, find
the moduli space of complex structures on $Y$, blow up to obtain a model
in which the boundary is a divisor with normal crossings, find the boundary
points with maximally unipotent monodromy, and sort them into equivalence
classes (as indicated in
section \ref{s:ambiguity}).  For one representative of each
equivalence class, find the canonical coordinates $z_j$
and the canonical gauge $\Omega(z)$
(these are unique if the integral monodromy
conjecture holds, otherwise make a choice) ,
calculate the three-point functions
\[\big\langle
\frac{\partial}{\partial z_j}
\frac{\partial}{\partial z_k}
\frac{\partial}{\partial z_\ell}
\big\rangle :=
\int_Y\Omega(z)\wedge\nabla_{z_j}\nabla_{z_k}\nabla_{z_\ell}\Omega(z)\]
in those coordinates and that gauge, and make a power series expansion
\[(2\pi i)^3 z_j z_k z_\ell \int_Y\Omega(z)\wedge\nabla_{z_j}
\nabla_{z_k}\nabla_{z_\ell}\Omega(z)
= c^{jk\ell} + \sum_{\eta\in H_2} c_\eta^{jk\ell} z^\eta, \]
where the leading term
$c^{jk\ell}$ should coincide with the intersection numbers on
the mirror partner, and
where we use a kind of multi-index notation for monomials $z^\eta$.

The coefficients
$c_\eta^{jk\ell}$ themselves are not quite the predictions for
``numbers of rational curves.''
One must take into account the ``multiple cover formula''  for the A-model
\cite{tftrc,Manin}\footnote{Although this multiple cover formula---first
postulated in \cite{CdGP}---can be derived using path integral arguments
 in conformal field theory \cite{tftrc} and has been mathematically proved
for the algebraic Gromov--Witten invariants \cite{Manin}, it is still not known
for the
symplectic Gromov--Witten invariants.}
and write the three-point functions in the form
\[(2\pi i)^3 z_j z_k z_\ell \int_Y\Omega(z)\wedge\nabla_{z_j}
\nabla_{z_k}\nabla_{z_\ell}\Omega(z)
= c^{jk\ell} + \sum_{\eta\in H_2} \eta_j\eta_k\eta_\ell\,\varphi_\eta \,
\frac{z^\eta}{1-z^\eta}. \]
The coefficients $\varphi_\eta$ are then the predicted number of rational
curves
in the homology class $\eta$.

\ifx\undefined\bysame
\newcommand{\bysame}{\leavevmode\hbox to3em{\hrulefill}\,}
\fi

\trailer

\end{document}